Analysis.

# Some explicit formulas for a sequence of secondary measures.

Roland Groux

*Lycée Polyvalent Rouvière, rue Sainte Claire Deville,
BP 1205. 83070 Toulon .Cedex. France.*

Email : roland.groux@orange.fr

**Abstract.**

We study here a sequence of secondary measures, so called because the set of secondary polynomials on a given term become orthogonal for the next measure. The main result is a formula making explicit the density of any term of the sequence, under some hypotheses. We give some applications and also derive an interpretation of the Fourier coefficients as multiple integrals.

## 1. Introduction and notations.

We consider a probability density function $x \mapsto \rho(x)$ on an interval $I$ bounded with $a$ and $b$. The Stieltjes transform of the measure of density $\rho$ is defined on $\mathbb{C} - I$ by the formula:

$$z \mapsto S_\rho(z) = \int_a^b \frac{\rho(t)dt}{z-t}. \ [1].$$

We note $n \mapsto P_n$ an Hilbert base of normalized polynomials for the classic inner product :

$(f,g) \mapsto <f/g>_\rho = \int_a^b f(t)g(t)\rho(t)dt$ on the associated Hilbert's space $L^2(I,\rho)$.

We call $Q_n(X) = \int_a^b \frac{P_n(t) - P_n(X)}{t-X}\rho(t)dt$ the secondary polynomial associated with $P_n$.

Let us recall the well known result below : [2]

If a positive measure on $I$ associated to a density function $\mu$ having Stietjes's transformation gived by the formula: $S_\mu(z) = z - c_1 - \dfrac{1}{S_\rho(z)}$, then secondary polynomials $Q_n$ form an orthogonal family for the inner product induced by $\mu$.

We habusively call $\mu$ the secondary measure associated with $\rho$. The moment of order 0 of this new measure is given by $d_0 = c_2 - (c_1)^2$.

If we normalize $\mu$, we introduce $\rho_1 = \dfrac{\mu}{d_0}$ and so we can continue the process with $\rho_1$. It thus appears a sequence of probability density functions : $n \mapsto \rho_n$ starting with $\rho_0 = \rho$ and such that for every *n* integer , $\rho_{n+1}$ is the 'normalized secondary measure' of $\rho_n$.

We adopt now the following notations :

$$c_1^n = \int_a^b x \rho_n(x)dx \ ; \ c_2^n = \int_a^b x^2 \rho_n(x)dx \ \text{ and } \ d_0^n = c_2^n - (c_1^n)^2 \ .$$

The Stieltjes transform of $\rho_n$ will be simply represented by $z \mapsto S_n(z)$.

## 2. A classical scheme of continuous fractions.

According with previous notations, we have for every $n$: $S_{n+1}(z) = \dfrac{1}{d_0^n}[z - c_1^n - \dfrac{1}{S_n(z)}]$.

This can also be written as : $S_n(z) = \dfrac{1}{z - c_1^n - d_0^n S_{n+1}}$. Thus we obtain a diagram of continuous fractions : $S_\rho(z) = S_0(z) = \cfrac{1}{z - c_1^0 - \cfrac{d_0^0}{z - c_1^1 - \cfrac{d_0^1}{z - c_1^2 - \cfrac{d_0^2}{z - \ldots \cfrac{\ldots}{z - c_1^{n-1} - \cfrac{d_0^{n-1}}{z - c_1^n - d_0^n S_{n+1}(z)}}}}}}$

According to classical results of the theory [3,4] we have the formula (2.1) :

- $S_\rho(z) = \dfrac{u_n(-d_0^n S_{n+1}(z)) + u_{n+1}}{v_n(-d_0^n S_{n+1}(z)) + v_{n+1}}$ , with the relations (2.2) $\begin{cases} u_{n+1} = (z - c_1^n)u_n - d_0^{n-1}u_{n-1} \\ v_{n+1} = (z - c_1^n)v_n - d_0^{n-1}v_{n-1} \end{cases}$

starting with : $\begin{cases} u_0 = 0 \\ u_1 = 1 \end{cases}$ ; $\begin{cases} v_0 = 1 \\ v_1 = z - c_1^0 \end{cases}$

- The determinant $\Delta_{n+1} = \begin{vmatrix} u_{n+1} & v_{n+1} \\ u_n & v_n \end{vmatrix}$ simplifies to : $\Delta_{n+1} = d_0^0 d_0^1 \ldots d_0^{n-1}$. (2.3).

We deduce immediatly that : $S_\rho(z) - \dfrac{u_{n+1}}{v_{n+1}} = \dfrac{(d_0^0 d_0^1 \ldots d_0^n) . S_{n+1}(z)}{v_{n+1}(v_{n+1} - d_0^n v_n S_{n+1}(z))}$ .

Let us study this difference in the neightborhood of infinity.

$S_{n+1}(z) = S_{\rho_{n+1}}(z)$ is equivalent to $\dfrac{1}{z}$. ( $\rho_{n+1}$ is a density of probability).

$v_{n+1}$ is obviously a polynomial in $z$ with leading term equal to $z^{n+1}$.

Thus, in the neightborhood of infinity : $S_\rho(z) - \dfrac{u_{n+1}}{v_{n+1}}$ is equivalent to $\dfrac{d_0^0 d_0^1 \ldots d_0^n}{z^{2n+3}}$

Noting that the degree of $u_{n+1}$ is $n$ and that of $v_{n+1}$ is $n+1$, we can conclude that the fraction

$F_n(z) = \dfrac{u_{n+1}}{v_{n+1}}$ is a Pade approximant for $S_\rho(z)$ of type $[n/n+1]$. [4]

According the classic theory, this fraction is simply equal to the quotient $\dfrac{Q_{n+1}(z)}{P_{n+1}(z)}$. So we have the proportionnality: $v_n = \lambda_n P_n(z)$ et $u_n = \lambda_n Q_n(z)$.

Thanks to the recurrence relation, the leading coefficient of $v_n$ is equal to 1, for every integer $n$, so we have $\lambda_n = \dfrac{1}{a_n}$, with $a_n$ leading coefficient of $P_n$.

From the formula above (2.1) : $S_\rho(z) = \dfrac{u_n(-d_0^n S_{n+1}(z)) + u_{n+1}}{v_n(-d_0^n S_{n+1}(z)) + v_{n+1}}$ and we can now deduce

$$S_n(z) = \frac{v_n S_0(x) - u_n}{d_0^{n-1}(v_{n-1} S_0(x) - u_{n-1})} = \frac{1}{d_0^{n-1}} \left(\frac{a_{n-1}}{a_n}\right) \frac{Q_n(z) - P_n(z) S_\rho(z)}{Q_{n-1}(z) - P_{n-1}(z) S_\rho(z)} \quad (2.4)$$

## 3. Calculation of coefficients $d_0^n$.

Another result of the theory of Stieltjes transform states that in the neightborhood of infinity :

$S_\rho(z) - \dfrac{Q_n(z)}{P_n(z)}$ is equivalent to $\dfrac{\gamma_n}{z^{2n+1}}$, with $\gamma_n = \dfrac{\gamma_0 a_0^2}{a_n^2} \times \dfrac{\|P_n\|^2}{\|P_0\|^2}$, expression in which $a_n$ is the leading coefficient of $P_n$. ( $\gamma_n = \dfrac{1}{a_n}\int_a^b t^n P_n(t)\rho(t)dt$ ). [6]

Comparing with the previously obtained equivalence, we deduce :

$$d_0^0 d_0^1 ........d_0^{n-1} = \frac{\gamma_0 a_0^2}{a_n^2} \times \frac{\|P_n\|^2}{\|P_0\|^2} \text{, and so : } d_0^n = \frac{a_n^2 \|P_{n+1}\|^2}{a_{n+1}^2 \|P_n\|^2} = \frac{a_n^2}{a_{n+1}^2}. \quad (3.1)$$

The formula (2.4) simplifies then to : (3.2) $S_n(z) = \dfrac{a_n}{a_{n-1}} \dfrac{Q_n(z) - P_n(z) S_\rho(z)}{Q_{n-1}(z) - P_{n-1}(z) S_\rho(z)}$

Recall now that in the classical recurrence relation to three terms, as written

(3.3) $xP_n(x) = t_n P_{n+1}(x) + s_n P_n(x) + t_{n-1} P_{n-1}(x)$, we have $t_{n-1} = \dfrac{a_{n-1} \|P_n\|^2}{a_n \|P_{n-1}\|^2} = \dfrac{a_{n-1}}{a_n}$, and so we

obtain (3.4) $S_n(z) = \dfrac{1}{t_{n-1}} \dfrac{Q_n(z) - P_n(z) S_\rho(z)}{Q_{n-1}(z) - P_{n-1}(z) S_\rho(z)}$

# 4. Explicitation of the density $\rho_n$.

We recall the inversion formula of Stieltjes Perron which allows to reconstruct the density from its Stieltjes transform under some hypotheses : $\rho(x) = \dfrac{1}{2i\pi}\lim\limits_{\varepsilon \to 0^+}(S_\rho(x-i\varepsilon) - S_\rho(x+i\varepsilon))$.

A simple case of applications appears when the density is a continuous function over a compact interval.

We suppose here the initial density $\rho$ satisfying the hypotheses of inversion and we suppose also the existence of $\varphi(x) = \lim\limits_{\varepsilon \to 0^+} S_\rho(x-i\varepsilon) + S_\rho(x+i\varepsilon)$ for $x$ all over the interval $I$.

This function is called the reducer of the measure of density $\rho$ and allows to explicit the secondary measure associated with $\rho$ as : $\mu(x) = \dfrac{\rho(x)}{\dfrac{\varphi^2(x)}{4} + \pi^2 \rho^2(x)}$. [5]

Under these hypotheses, the density $\rho_n$ can be made explicit thanks Stieltjes Perron and the formula (3.4) above : $S_n(z) = \dfrac{1}{t_{n-1}} \dfrac{Q_n(z) - P_n(z) S_\rho(z)}{Q_{n-1}(z) - P_{n-1}(z) S_\rho(z)}$.

An elementary calculation then leads to :

$$\rho_n(x) = \dfrac{1}{t_{n-1}} \times \dfrac{\rho(x)(P_{n-1}(x)Q_n(x) - P_n(x)Q_{n-1}(x))}{(P_{n-1}(x)\dfrac{\varphi(x)}{2} - Q_{n-1}(x))^2 + \pi^2 \rho^2(x) P_{n-1}^2(x)} \quad (4.1)$$

Now recall the formula (2.2) : $\Delta_{n+1} = \begin{vmatrix} u_{n+1}(x) & v_{n+1}(x) \\ u_n(x) & v_n(x) \end{vmatrix} = d_0^0 . d_0^1 ...... d_0^{n-1}$

From $v_n = \lambda_n P_n(z)$ ; $u_n = \lambda_n Q_n(z)$ ; $\lambda_n = \dfrac{1}{a_n}$ ; $d_0^n = \dfrac{a_n^2}{a_{n+1}^2}$ we easily deduce :

$Q_{n+1}(x) P_n(x) - P_{n+1}(x) Q_n(x) = \dfrac{a_{n+1}}{a_n} = \dfrac{1}{t_n}$. So, the formula (4.1) simplifies to :

$$(4.2) \quad \boxed{\rho_n(x) = \dfrac{1}{(t_{n-1})^2} \times \dfrac{\rho(x)}{(P_{n-1}(x)\dfrac{\varphi(x)}{2} - Q_{n-1}(x))^2 + \pi^2 \rho^2(x) P_{n-1}^2(x)}}$$

Since all these functions are probability densities, we deduce the value of the following integral as :

(4.3) $$\int_a^b \frac{\rho(x)}{(P_{n-1}(x)\frac{\varphi(x)}{2} - Q_{n-1}(x))^2 + \pi^2 \rho^2(x) P_{n-1}^2(x)} dx = (t_{n-1})^2$$

Some examples. [5,6,7]

- For the uniform Lebesgue measure over [0,1] :

$$\int_0^1 \frac{dx}{[P_n(x)\ln(\frac{x}{1-x}) - Q_n(x)]^2 + \pi^2 P_n^2(x)} = \frac{(n+1)^2}{4(2n+1)(2n+3)}.$$

$P_n(x)$ is the shifted (and normalized) Legendre polynomial of order $n$.

- For $\rho(x) = e^{-x}$ over $[0, +\infty[$

$$\int_0^{+\infty} \frac{e^{-x} dx}{[P_n(x) e^{-x} \text{Ei}(x) - Q_n(x)]^2 + \pi^2 e^{-2x} P_n^2(x)} = (n+1)^2$$

$P_n(x)$ is the Laguerre polynomial of order $n$ and Ei the exponential integral.

- For the Gaussian measure $\rho(x) = \dfrac{e^{-\frac{x^2}{2}}}{\sqrt{2\pi}}$ over the real axis.

$$\int_{-\infty}^{+\infty} \frac{\rho(x) dx}{[\sqrt{\frac{\pi}{2}} e^{-\frac{x^2}{2}} \text{erfi}(\frac{x}{\sqrt{2}}) P_n(x) - Q_n(x)]^2 + \pi^2 \rho^2(x) P_n^2(x)} = n+1$$

$P_n(x)$ is the Hermite polynomial of order $n$ and erfi the imaginary error function.

- For the Tchebychev measure of the second kind ( $\rho(x) = \dfrac{2}{\pi}\sqrt{1-x^2}$ over [-1,1])

In this case the sequence of normalized secondary measures is constant. So we find the classical relation connecting the Tchebychev polynomials :

$$P_n^2(x) - P_{n-1}(x) P_{n+1}(x) = 1$$

## 5. An interpretation of the Fourier coefficients.

We recall the followin result : [5]

The operator $f(x) \mapsto g(x) = \int_a^b \frac{f(t) - f(x)}{t - x} \rho(t) dt$ creating secondary polynomials extends to a continuous linear map $T_\rho$ linking the space $L^2(I, \rho)$ to the Hilbert'space $L^2(I, \mu)$.

Its restriction to the hyperplane $H_\rho$ of the function orthogonal with $P_0 = 1$ constitutes an isometric function for both norms respectively. So, for any couple $(f, g)$ of elements of $L^2(I, \rho)$, we have:

(5.1) $\boxed{<f/g>_\rho - <f/1>_\rho \times <g/1>_\rho = <T_\rho(f)/T_\rho(g)>_\mu}$

Because of the normalization performed in the sequence of secondary measures, we lose the isometric character of the asociated transforms.

If we note $H_n$ the space of elements of $L^2(I, \rho_n)$ orthogonal with $P_0 = 1$,

$\sqrt{d_0^n} \times T_{\rho_n}$ is now an isometric map linking $(H_n, \rho_n)$ to $(L^2(I, \rho_{n+1})$

So, for any couple $(f, g)$ of $L^2(I, \rho_n)$, we can write :

(5.2) $\boxed{<f/g>_{\rho_n} - <f/1>_{\rho_n} \times <g/1>_{\rho_n} = d_0^n <T_{\rho_n}(f)/T_{\rho_n}(g)>_{\rho_{n+1}}}$

We will now generalize this formula by composing the transforms.

We introduce $F_n = T_{\rho_n} \circ T_{\rho_{n-1}} \circ ..... \circ T_{\rho_1} \circ T_{\rho_0}$. ( $\rho_0 = \rho$ and so $F_0 = T_\rho$ )

An elementary recurrence led to :

(5.3) $F_n(f)(t_{n+1}) = \int_a^b \int_a^b ..... \int_a^b (\sum_{k=0}^{k=n+1} \frac{f(t_k)}{\prod_{j \neq k}(t_k - t_j)}) \rho_0(t_0)...\rho_n(t_n) dt_0....dt_n$.

(multiple integral of order $n+1$)

➢ Let be a couple $(f, g)$ of elements of $L^2(I, \rho)$ with $g \perp P_0$.
we deduce then from (5.2) : $<f/g>_{\rho_0} = <T_{\rho_0}(f)/T_{\rho_0}(g)>_\mu = d_0^0 <F_0(f)/F_0(g)>_{\rho_1}$

➢ Assume further $g$ orthogonal with $P_1$, still in $L^2(I, \rho)$. We deduce that $F_0(g)$ is orthogonal with $P_0 = 1$ in $L^2(I, \rho_1)$, and so we get : $<f/g>_{\rho_0} = d_0^0 d_0^1 <F_1(f)/F_1(g)>_{\rho_2}$

➢ More generally, if $g$ is orthogonal at every term of $(P_0, P_1, ...., P_n)$ in $L^2(I, \rho)$.

We have the relation (5.4) : $\boxed{<f/g>_{\rho_0} = d_0^0 d_0^1 .....d_0^n <F_n(f)/F_n(g)>_{\rho_{n+1}}}$

This formula above will allow us to make explicit the Fourier coefficients of a given function as multiple integral.

Because $P_n \perp P_0, P_1, ..... P_{n-1}$ we deduce directly from (5.4):

(5.5) $C_n(f) = <f/P_n>_{\rho_0} = d_0^0......d_0^{n-1} <F_{n-1}(f)/F_{n-1}(P_n)>_{\rho_n}$

Recall now a precedent result (3.1): $d_0^0.....d_0^{n-1} = \frac{\gamma_0 a_0^2}{a_n^2} \times \frac{\|P_n\|^2}{\|P_0\|^2} = \frac{1}{a_n^2}$

We also easilty result $F_{n-1}(P_n) = a_n$, leading coefficient of $P_n$, because every function $\rho_n$ is a density of probability.

So the formula (5.5) simplifies to :

(5.6) $$\boxed{C_n(f) = \frac{1}{a_n} \int_a^b \int_a^b .... \int_a^b (\sum_{k=0}^{k=n} \frac{f(t_k)}{\prod_{j \neq k}(t_k - t_j)}) \rho_0(t_0) .....\rho_n(t_n) dt_0 .... dt_n}$$

We conclude this paragraph by a particular application of this formula.

We consider here the function : $x \mapsto f(x) = \frac{1}{x+a}$, with $a \notin I$

$f$ is an eigenvector for every operator $T_\rho$. Specifically, an elementary calculation gives :

$T_\rho(f(x)) = \gamma_a f(x)$ with $\gamma_a = -\int_I \frac{\rho(t)dt}{t+a} = S_\rho(-a)$

By elementary composition, we get : $F_{n-1}(f) = \gamma_0 \times \gamma_1 \times ..... \times \gamma_{n-1} \times f$. So, from (5.6) :

$C_n(f) = \frac{1}{a_n} \times \gamma_0 \times \gamma_1 \times ..... \times \gamma_{n-1} \times \int_I \frac{\rho_n(t)}{t+a} dt = -\frac{1}{a_n} \times \gamma_0 \times \gamma_1 \times ..... \times \gamma_{n-1} \times \gamma_n$

This can be written as: $\boxed{C_n(f) = -\frac{1}{a_n} \times \prod_{k=0}^{k=n} S_{\rho_k}(-a)}$

And more explicitly : $\int_I \frac{P_n(x)\rho(x)dx}{x+a} = -\frac{1}{a_n} \times \prod_{k=0}^{k=n} S_{\rho_k}(-a)$

By simple quotient : $S_{\rho_n}(-a) = \frac{a_n}{a_{n-1}} \times \frac{C_n(f)}{C_{n-1}(f)}$

So we have for $z \notin I$ : $\boxed{\int_I \frac{\rho_n(t)}{z-t} dt = \frac{a_n}{a_{n-1}} \times \frac{\int_I \frac{P_n(t)\rho(t)dt}{t-z}}{\int_I \frac{P_{n-1}(t)\rho(t)dt}{t-z}}}$ (5.7)

Apply this formula to the tchebychev measure of the second kind over the interval $[0, 1]$.

We have $\rho(x) = \dfrac{8\sqrt{x(1-x)}}{\pi}$, $S_{\rho_t}(-a) = -8a - 4 + 8\sqrt{a^2 + a}$, and the leading coefficient of $P_n$ is $a_n = 4^n$. As in this case the sequence of normalized secondary measures is constant, we deduce from

$$(5.7): -8a - 4 + 8\sqrt{a^2 + a} = 4 \times \dfrac{\int_0^1 \dfrac{P_n(t)\rho_t(t)dt}{t+a}}{\int_0^1 \dfrac{P_{n-1}(t)\rho_t(t)dt}{t+a}}$$

$C_n(f) = q^n C_0(f)$, with $q = 2\sqrt{a^2 + a} - 2a - 1$ and $C_0(f) = \int_0^1 \dfrac{\rho_t(x)dx}{x+a} = -4q$

So we have in $L^2([0,1], \rho)$ $\quad \dfrac{1}{x+a} = -4 \sum_{n=0}^{n=\infty} q^{n+1} P_n(x)$.

Note then $q = 2\sqrt{a^2 + a} - 2a - 1 \Rightarrow a = -\dfrac{(q+1)^2}{4q}$

So we get : $\boxed{\dfrac{t}{(t+1)^2 - 4tx} = \sum_{n=0}^{n=\infty} t^{n+1} P_n(x)}$

One recognizes here the classical generating function for the Tchebychev polynomials.

## 6. Complements and other formulas.

- First moment of $\rho_n$.

From the formulas (2.2) we get : $v_{n+1}(x) = (x - c_1^n)v_n(x) - d_0^{n-1} v_{n-1}(x)$.

Thanks to the proportionality $v_n = \lambda_n P_n$, this can be written :

$$xP_n(x) = \dfrac{\lambda_{n+1}}{\lambda_n} P_{n+1}(x) + c_1^n P_n(x) + \dfrac{d_0^{n-1}\lambda_{n-1}}{\lambda_n} P_{n-1}(x)$$

If we recall the classical three-terms relation (3.3) :

$xP_n(x) = t_n P_{n+1}(x) + s_n P_n(x) + t_{n-1} P_{n-1}(x)$, we deduce immediatly the equalities :

$c_1^n = s_n$ ; $t_n = \dfrac{\lambda_{n+1}}{\lambda_n}$ and $t_{n-1} = \dfrac{d_0^{n-1}\lambda_{n-1}}{\lambda_n}$

Note that the last two give back : $d_0^n = (t_n)^2$.

Using (4.2), the relation $c_1^n = s_n$ is reflected in integral form by

$$\boxed{\int_a^b \dfrac{x\rho(x)}{(P_{n-1}(x)\dfrac{\varphi(x)}{2} - Q_{n-1}(x))^2 + \pi^2 \rho^2(x) P_{n-1}^2(x)} dx = s_n (t_{n-1})^2} \quad (6.1)$$

- Reducer of $\rho_n$.

Using formula (3.4) we can make explicit the reducer of $\rho_n$, defined by :

$\varphi_n(x) = \lim_{\varepsilon \to 0^+} S_n(x - i\varepsilon) + S_n(x + i\varepsilon)$. After some elementary calculations and thanks to the inversion formula of Stieltjes Perron, we get :

$$\varphi_n(x) = \frac{2}{(t_{n-1})^2} \times \frac{(P_n(x)\frac{\varphi(x)}{2} - Q_n(x)) \times (P_{n-1}(x)\frac{\varphi(x)}{2} - Q_{n-1}(x)) + \pi^2 \rho^2(x) P_n(x) P_{n-1}(x)}{(P_{n-1}(x)\frac{\varphi(x)}{2} - Q_{n-1}(x))^2 + \pi^2 \rho^2(x) P_{n-1}^2(x)} \quad (6.2)$$

- **Orthogonal polynomials associated with $\rho_n$.**

From the coupling formula $S_\mu(z) = z - c_1 - \frac{1}{S_\rho(z)}$ we can easily express the set of primary $(U_n)$ and secondary $(V_n)$ orthogonal polynomials of $\mu$ through :

$$\begin{cases} U_n(x) = Q_{n+1}(x) \\ V_n(x) = (x - c_1)Q_{n+1}(x) - P_{n+1}(x) \end{cases} \quad (6.3)$$

From the isometric character of $T_\rho$, we have for $n \geq 1$ $\|Q_n\|_\mu = \|P_n\|_\rho = 1$.

When we normalize $\mu$ to $\rho_1 = \frac{\mu}{d_0}$, we can introduce $P_n^1(x) = t_0 U_n(x)$ to keep the norm equal to 1. ( Because $T_{\rho_1} = \frac{1}{d_0} T_\mu$ and $d_0 = (t_0)^2$ )

The same way we change $V_n$ to $Q_n^1(x) = \frac{V_n(x)}{t_0}$ .

So we get $Q_n^1 = T_{\rho_1}(P_n^1)$ because $T_{\rho_1} = \frac{1}{d_0} T_\mu$ and $T_\mu(U_n) = V_n$.

The relations (6.3) translate in to : $\begin{cases} P_n^1(x) = t_0 Q_{n+1}(x) \\ Q_n^1(x) = \frac{1}{t_0}[(x - c_1)Q_{n+1}(x) - P_{n+1}(x)] \end{cases} \quad (6.4)$

More genrally if we note $n \mapsto (P_n^k(x), Q_n^k(x))$ the set of polynomials associated with $\rho_k$, we can write in matrix form : :$\begin{pmatrix} P_n^{k+1}(x) \\ Q_n^{k+1}(x) \end{pmatrix} = \frac{1}{t_k}\begin{pmatrix} 0 & (t_k)^2 \\ -1 & x - c_1^k \end{pmatrix}\begin{pmatrix} P_{n+1}^k(x) \\ Q_{n+1}^k(x) \end{pmatrix}$ (6.5)

Now introduce $M_k(x) = \dfrac{1}{t_k}\begin{pmatrix} 0 & t_k^2 \\ -1 & x - s_k \end{pmatrix}$ and $\Pi_k(x) = M_k(x).M_{k-1}(x).....M_0(x)$

We easily deduce from (6.5) : $\begin{pmatrix} P_n^{k+1}(x) \\ Q_n^{k+1}(x) \end{pmatrix} = \Pi_k(x) \begin{pmatrix} P_{n+1+k}^0(x) \\ Q_{n+1+k}^0(x) \end{pmatrix}$ (6.6)

If we note $\Pi_k(x) = \begin{pmatrix} A_{k+1}(x) & B_{k+1}(x) \\ C_{k+1}(x) & D_{k+1}(x) \end{pmatrix}$, we obtain the relations :

$$\begin{cases} A_{k+2}(x) = t_{k+1} C_{k+1}(x) & B_{k+2}(x) = t_{k+1} D_{k+1}(x) \\ C_{k+2}(x) = \dfrac{-1}{t_{k+1}} A_{k+1}(x) + \dfrac{(x - s_{k+1})}{t_{k+1}} C_{k+1}(x) & D_{k+2}(x) = \dfrac{-1}{t_{k+1}} B_{k+1}(x) + \dfrac{(x - s_{k+1})}{t_{k+1}} D_{k+1}(x) \end{cases}$$

It appears clear that $C_n$ and $D_n$ satisfy the classical three-terms relation for the initial measure : $xP_n(x) = t_n P_{n+1}(x) + s_n P_n(x) + t_{n-1} P_{n-1}(x)$

$\forall k \geq 1 \quad xC_k(x) = t_k C_{k+1}(x) + s_k C_k(x) + t_{k-1} C_{k-1}(x)$

$\forall k \geq 1 \quad xD_k(x) = t_k D_{k+1}(x) + s_k D_k(x) + t_{k-1} D_{k-1}(x)$

According to the classical theory, we can write : $\begin{cases} C_k(x) = aP_k(x) + bQ_k(x) \\ D_k(x) = cP_k(x) + dQ_k(x) \end{cases}$

The components $a, b, c, d$ are determined using the initial conditions.

The sequence starts with $\begin{cases} P_0^0(x) = 1; \quad P_1^0(x) = \dfrac{x - s_0}{t_0} \\ Q_0^0(x) = 0; \quad Q_1^0(x) = \dfrac{1}{t_0} \end{cases}$

By (6.4) we have $P_2^0(x) = \dfrac{(x - s_0)(x - s_1) - (t_0)^2}{t_0 t_1}$; $Q_2^0(x) = \dfrac{x - s_1}{t_0 t_1}$.

A direct calculation with matrices gives us :

$C_1(x) = -Q_1^0(x); \quad C_2(x) = -Q_2^0(x); \quad D_1(x) = P_1^0(x); \quad D_2(x) = P_2^0(x)$

So we get : $a = 0$ ; $b = -1$ ; $c = 1$ ; $d = 0$ and consequently :

For $k \geq 1 \quad C_k(x) = -Q_k^0(x) = -Q_k(x)$ and $D_k(x) = P_k^0(x) = P_k(x)$

$\Pi_k(x) = \begin{pmatrix} -t_k Q_k(x) & t_k P_k(x) \\ -Q_{k+1}(x) & P_{k+1}(x) \end{pmatrix}$ ; $\boxed{\begin{cases} P_n^{k+1}(x) = t_k [P_k(x) Q_{n+k+1}(x) - Q_k(x) P_{n+k+1}(x)] \\ Q_n^{k+1}(x) = P_{k+1}(x) Q_{n+k+1}(x) - P_{n+k+1}(x) Q_{k+1}(x) \end{cases}}$ (6.7)